\documentclass[10pt]{amsart}
\usepackage{amsfonts}
\usepackage{amssymb}
\usepackage{amsthm}
\usepackage{latexsym}
\usepackage{amsmath}
\usepackage{amssymb}

\newtheorem{thm}{Theorem}
\newtheorem{prop}[thm]{Proposition}
\newtheorem{claim}[thm]{Claim}
\newtheorem{obs2}[thm]{Observation}
\newtheorem{cor}[thm]{Corollary}
\newtheorem{lemma}[thm]{Lemma}
\newtheorem{defin}{Definition}
\newenvironment{dem}{\begin{proof}[Proof]}{\end{proof}}

\usepackage{amsfonts}
\newcommand{\field}[1]{\mathbb{#1}}
\newcommand{\Q}{\field{Q}}

\newcommand{\R}{\field{R}}
\newcommand{\N}{\field{N}}
\newcommand{\Z}{\field{Z}}

\begin{document}

\title[Normalizer of $\Gamma_0(N)$]{The group structure of the normalizer of $\Gamma_0(N)$}

\author{Francesc Bars }

%\thanks{Supported by BFM2006-????.}

\date{First version: December 22th, 2006}
\thanks{MSC: 20H05(19B37,11G18)}
\maketitle

\begin{center}
\begin{small}
\begin{abstract}
We determine the group structure of the normalizer of
 $\Gamma_0(N)$ in $SL_2(\R)$ modulo $\Gamma_0(N)$. These results correct the Atkin-Lehner statement
 \cite[Theorem 8]{AL}.
\end{abstract}
\end{small}
\end{center}

\section{Introduction}
The modular curves $X_0(N)$ contain deep arithmetical information.
These curves are the Riemann surfaces obtained by completting with
the cusps the upper half plane modulo the
modular subgroup $$\Gamma_0(N)=\{\left(\begin{array}{cc} a&b\\
Nc&d\\
\end{array}\right)\in SL_2(\mathbb{Z})|c\in\Z\}.$$
It is clear that the elements in the normalizer of $\Gamma_0(N)$
in $SL_2(\R)$ induce automorphisms of $X_0(N)$ and moreover one
obtains in that way all automorphisms of $X_0(N)$ for $N\neq 37$
and $63$ \cite{KM}. This is one reason coming from the modular
world that shows the interest in computing the group structure of
this normalizer modulo $\Gamma_0(N)$.

Morris Newman obtains a result for this normalizer in terms of
matrices \cite{N},\cite{N2}, see also the work of Atkin-Lehner and
Newman \cite{NL}. Moreover, Atkin-Lehner state without proof the
group structure of this normalizer modulo $\Gamma_0(N)$
\cite[Theorem 8]{AL}. In this paper we correct this statement and
we obtain the right structure of the normalizer modulo
$\Gamma_0(N)$. The results are a generalization of some results
obtained in \cite{Ba}.

\section{The Normalizer of $\Gamma_0(N)$ in $SL_2(\R)$}

Denote by $\operatorname{Norm}(\Gamma_0(N)$ the normalizer of
$\Gamma_0(N)$ in $SL_2(\R)$.

\begin{thm}[Newman]\label{teo1} Let $N=\sigma^2 q$ with $\sigma,q\in\N$ and $q$ square-free.
Let $\epsilon$ be the $\operatorname{gcd}$ of all integers of the
form $a-d$ where $a,d$ are integers such that $\left(\begin{array}{cc} a&b\\
Nc&d\\
\end{array}\right)\in\Gamma_0(N)$. Denote by $v:=v(N):=\operatorname{gcd}(\sigma,\epsilon)$.
Then $M\in \operatorname{Norm}(\Gamma_0(N))$ if and only if $M$ is
of the form
$$\sqrt{\delta}\left(\begin{array}{cc} r\Delta&\frac{u}{v\delta\Delta}\\
\frac{sN}{v\delta\Delta}&l\Delta\\
\end{array}\right)$$
with $r,u,s,l\in\Z$ and $\delta|q$, $\Delta|\frac{\sigma}{v}$.
Moreover $v=2^{\mu}3^{w}$ with $\mu=min(3,[\frac{1}{2}v_2(N)])$
and $w=min(1,[\frac{1}{2}v_3(N)])$ where $v_{p_i}(N)$ is the
valuation at the prime $p_i$ of the integer $N$.
\end{thm}
This theorem is proved by Morris Newman in \cite{N} \cite{N2}, see
also \cite[p.12-14]{Ba}.

Observe that if $\operatorname{gcd}(\delta\Delta,6)=1$ we have
$\operatorname{gcd}(\delta\Delta^2,\frac{N}{\delta\Delta^2})=1$
because the determinant is one .

\section{The group structure of $\operatorname{Norm}(\Gamma_0(N))/\Gamma_0(N)$}

In this section we obtain some partial results on the group
structure of $\operatorname{Norm}(\Gamma_0(N))$. Let us first
introduce some particular elements of $SL_2(\R)$.
\begin{defin}  Let $N$ be fixed. For every divisor $m'$ of $N$ with
$\operatorname{gcd}(m',N/m')=1$ the Atkin-Lehner involution
$w_{m'}$ is defined as follows,
$$w_{m'}=\frac{1}{\sqrt{m'}}\left(\begin{array}{cc} m'a&b\\
Nc&m'd\\
\end{array}\right)\in SL_2(\R)$$
with $a,b,c,d\in\Z$.
\end{defin}

Denote by $S_{v'}=\left(\begin{array}{cc} 1&\frac{1}{v'}\\
0&1\\
\end{array}\right)$ with $v'\in\N\setminus\{0\}$. Atkin-Lehner claimed in \cite{AL} the
following:
\begin{claim}[Atkin-Lehner]\cite[Theorem 8]{AL}\label{AL}
The quotient $\operatorname{Norm}(\Gamma_{0}(N))/\Gamma_0(N)$ is
the direct product of the following groups:
\begin{enumerate}
\item \{$ w_{q^{\upsilon_q(N)}}$\} for every prime $q$, $q\ge5$
$q\mid N$. \item
\begin{enumerate}
      \item If $\upsilon_{3}(N)=0$, \{$1$\}
      \item If $\upsilon_{3}(N)=1$, \{$w_{3}$\}
      \item If $\upsilon_{3}(N)=2$, \{$w_{9},S_{3}$\}; satisfying $w_{9}^2=S_{3}^3=(w_{9}S_{3})^3=1$ (factor of order 12)
      \item If $\upsilon_{3}(N)\ge3$; \{$w_{3^{\upsilon_{3}(N)}},S_{3}$\}; where $w_{3^{\upsilon_{3}(N)}}^2=S_{3}^3=1$ and $w_{3^{\upsilon_{3}(N)}}S_{3}w_{3^{\upsilon_{3}(N)}}$
      commute
      with $S_{3}$ (factor group with 18 elements)
      \end{enumerate}
\item Let be $\lambda=\upsilon_{2}(N)$ and
$\mu=\operatorname{min}(3,[\frac{\lambda}{2}])$ and denote by
$\upsilon''=2^{\mu}$ the we have:
      \begin{enumerate}
      \item If $\lambda=0$ ; \{$1$\}
      \item If $\lambda=1$; \{$w_{2}$\}
      \item If $\lambda=2\mu$; \{$w_{2^{\upsilon_{2}(N)}},S_{\upsilon''}$\}
      with the relations
      $w_{2^{\upsilon_{2}(N)}}^2=S_{\upsilon''}^{\upsilon''}=(w_{2^{\upsilon_{2}(N)}}S_{\upsilon''})^3=1$,
      where they have orders 6,24, and 96 for $\upsilon=2,4,8$ respectively. (One needs to warn that for $v=8$ the relations
      do not define totally this factor group).
      \item If $\lambda> 2\mu$; \{ $w_{2^{\upsilon_{2}(N)}},S_{\upsilon''}$\}; $w_{2^{\upsilon_{2}(N)}}^2=S_{\upsilon''}^{\upsilon''}=1.$
       Moreover, $S_{\upsilon''}$ commutes with $w_{2^{\upsilon_{2}(N)}}S_{\upsilon''}w_{2^{\upsilon_{2}(N)}}$
       (factor group of order $2 {\upsilon''}^2$).
      \end{enumerate}
\end{enumerate}
\end{claim}

Let us give some partial results first.
\begin{prop}\label{propv1} Suppose that $v(N)=1$ (thus $4\nmid N$ and $9\nmid
N$). Then the Atkin-Lehner involutions generate
$\operatorname{Norm}(\Gamma_0(N)/\Gamma_0(N)$ and the group
structure is
$$\cong \prod_{i=1}^{\pi(N)}\Z/2\Z$$
where $\pi(N)$ is the number of prime numbers $\leq N$.
\end{prop}
\begin{dem} This is classically known. We recall only that $w_{m m'}=w_{m}w_{m'}$ for $(m,m')=1$ and
easily $w_mw_{m'}=w_{m'}w_m$; the the result follows by a
straightforward computation from Theorem \ref{teo1}, see also
\cite[p.14]{Ba}.
\end{dem}
When $v(N)>1$ it is clear that some element $S_{v'}$ appears in
the group structure of
$\operatorname{Norm}(\Gamma_0(N))/\Gamma_0(N)$ from Theorem
\ref{teo1}.
\begin{lemma}\label{lema2} If $4|N$ the involution $S_2\in \operatorname{Norm}(\Gamma_0(N))$
commutes with the Atkin-Lehner involutions $w_{m}$ with
$\operatorname{gcd}(m,2)=1$ and with the other $S_{v'}$.
\end{lemma}
\begin{dem} By the hypothesis the following matrix
belongs to $\Gamma_0(N)$
$$w_{m}S_2w_{m}S_2=\left(\begin{array}{cc} \frac{2mk^2+2Nt+mkNt}{2m}&\frac{(2+2m)(2m+2mk+Nt)}{4m}\\
\frac{Nt(2m+2mk+Nt)}{2m}&m+Nt+\frac{Nt}{m}+\frac{kNt}{2}+\frac{Nt^2}{4m}\\
\end{array}\right).$$
\end{dem}

\begin{prop}\label{propv2} Let $N=2^{v_2(N)}\prod_ip_i^{n_i}$, with $p_i$
different odd primes and assume that $v_2(N)\leq 3$, $v_3(N)\leq
1$. Then Atkin-Lehner's Claim \ref{AL} is true.
\end{prop}
For the proof we need two lemmas.

\begin{lemma} \label{nno:teo} Let $\tilde{u}\in \operatorname{Norm}(\Gamma_{0}(N))$ and write it as:
$$\tilde{u}=\frac{1}{\sqrt{\delta\Delta^2}}\left(\begin{array}{cc}
\Delta^2\delta r&\frac{u}{2}\\
\frac{sN}{2}&l\Delta^2\delta\\
\end{array}\right),$$
following the notation of Theorem \ref{teo1}.
 Then:\newline

$$w_{\Delta^2\delta}\tilde{u}=\left(\begin{array}{cc}
r'&\frac{u'}{2}\\
\frac{s'N}{2}&v'\\
\end{array}\right),\ if\ \operatorname{gcd}(\delta,2)=1,$$
$$w_{\Delta^2\frac{\delta}{2}}\tilde{u}=\frac{1}{\sqrt{2}}\left(\begin{array}{cc}
2r''&\frac{u''}{2}\\
\frac{s''N}{2}&2v''\\
\end{array}\right),\ if\ \operatorname{gcd}(\delta,2)=2.$$
\end{lemma}
\begin{dem}
This is an easy calculation. \end{dem}

We study now the different elements of the type
$$a(r',u',s',v')=\left(\begin{array}{cc}
r'&\frac{u'}{2}\\
\frac{s'N}{2}&v'\\
\end{array}\right),$$
$$b(r'',u'',s'',v'')=\frac{1}{\sqrt{2}}\left(\begin{array}{cc}
2r''&\frac{u''}{2}\\
\frac{s''N}{2}&2v''\\
\end{array}\right).$$
Observe that $b(,,,)$ only appears when $N\equiv0(mod\ 8)$.
\begin{lemma} \label{nna:teo} For $N\equiv4(mod\ 8)$ all
the elements of the normalizer of type $a(r',u',s',v')$ belong to
the order six group
$\{S_{2},w_{4}|S_{2}^2=w_{4}^2=(w_{4}S_{2})^3=1\}$.
\end{lemma}
\begin{dem} Straightforward from the equalities:
$$a(r',u',s',v')\in\Gamma_{0}(N)\Leftrightarrow s'\equiv u'\equiv0(mod\ 2)$$
$$a(r',u',s',v')S_{2}\in\Gamma_{0}(N)\Leftrightarrow r'\equiv v'\equiv\ u'\equiv1\ s'\equiv0(mod\ 2)$$
$$a(r',u',s',v')w_{4}\in\Gamma_{0}(N)\Leftrightarrow r'\equiv v'\equiv0\ u'\equiv s'\equiv1(mod\ 2)$$
$$a(r',u',s',v')w_{4}S_{2}\in\Gamma_{0}(N)\Leftrightarrow r'\equiv u'\equiv s'\equiv1\ v'\equiv0(mod\ 2)$$
$$a(r',u',s',v')S_{2}w_{4}\in\Gamma_{0}(N)\Leftrightarrow v'\equiv u'\equiv s'\equiv1\ r'\equiv0(mod\ 2)$$
$$a(r',u',s',v')S_{2}w_{4}S_{2}\in\Gamma_{0}(N)\Leftrightarrow r'\equiv v'\equiv s'\equiv1\ u'\equiv0(mod\ 2)$$

\end{dem}
\begin{lemma} \label{nne:teo} Let $N$ be a positive integer with $\upsilon_{2}(N)=3$. Then all the elements
of the form $a(r',u',s',v')$ and $b(r'',u'',s'',v'')$ correspond
to some element of the following group of 8 elements
$$\{ S_{2},w_{8}|S_{2}^2=w_{8}^2=1,S_{2}w_{8}S_{2}w_{8}=w_{8}S_{2}w_{8}S_{2}\}$$
\end{lemma}
\begin{dem} If follows from the equalities:
$$a(r',u',s',v')\in\Gamma_{0}(N)\Leftrightarrow r'\equiv v'\equiv1,u'\equiv s'\equiv0(mod\ 2)$$
$$a(r',u',s',v')S_{2}\in\Gamma_{0}(N)\Leftrightarrow r'\equiv v'\equiv u'\equiv1,s'\equiv0(mod\ 2)$$
$$a(r',u',s',v')w_{8}S_{2}w_{8}\in\Gamma_{0}(N)\Leftrightarrow r'\equiv v'\equiv s'\equiv1,u'\equiv0(mod\ 2)$$
$$a(r',u',s',v')S_{2}w_{8}S_{2}w_{8}\in\Gamma_{0}(N)\Leftrightarrow r'\equiv v'\equiv s'\equiv v'\equiv1(mod\ 2)$$
$$b(r'',u'',s'',v'')w_{8}\in\Gamma_{0}(N)\Leftrightarrow r''\equiv v''\equiv0,u''\equiv s''\equiv 1(mod\ 2)$$
$$b(r'',u'',s'',v'')S_{2}w_{8}S_{2}\in\Gamma_{0}(N)\Leftrightarrow r''\equiv v''\equiv u''\equiv s''\equiv1(mod\ 2)$$
$$b(r'',u'',s'',v'')S_{2}w_{8}\in\Gamma_{0}(N)\Leftrightarrow r''\equiv0,u''\equiv s''\equiv v''\equiv1(mod\ 2)$$
$$b(r'',u'',s'',v'')w_{8}S_{2}\in\Gamma_{0}(N)\Leftrightarrow v''\equiv0,u''\equiv s''\equiv r''\equiv1(mod\ 2)$$

\end{dem}
We can now proof Proposition \ref{propv2}].
\begin{dem}[ of Proposition \ref{propv2}]
Let $N=2^{\upsilon_{2}(N)}\prod_{i}p_{i}^{n_i}$, with $p_{i}$
different primes and assume that $9\nmid N$. If
$\upsilon_{2}(N)\le1$ we are done by proposition \ref{propv1}.
Suppose $\upsilon_{2}(N)=2$ and  let $\tilde{u}\in
\operatorname{Norm}(\Gamma_{0}(N))$. By lemmas \ref{nno:teo} and
\ref{nna:teo}, $w_{\delta}\tilde{u}=\alpha$, $\alpha\in\{
S_{2},w_{4}|S_{2}^2=w_{4}^2=(w_{4}S_{2})^3=1$ and it follows that
$\tilde{u}=w_{\delta}\alpha$. Since $w_{\delta}$ ($(\delta,2)=1$)
commutes with $S_{2}$ and the Atkin-Lehner involutions commute one
to each other, we are already done. In the situation $8||N$ the
proof is exactly the same but using lemmas \ref{nno:teo} and
\ref{nne:teo} instead.
\end{dem}

\section{Counterexamples to Claim \ref{AL}.}

In the above section we have seen that Atkin-Lehner's claim is
true if $v(N)\leq2$ i.e. for $v_2(N)\leq 3$ and $v_3(N)\leq 1$.
Now we obtain counterexamples when $v_2(N)$ and/or $v_3(N)$ are
bigger.

\begin{lemma} Claim \ref{AL} for $N=48$ is wrong.
\end{lemma}
\begin{dem} We know by Ogg \cite{O} that $X_0(48)$ is an
hyperelliptic modular curve with hyperelliptic involution not of
Atkin-Lehner type. The hyperelliptic involution always belongs to
the center of the automorphism group. We know by \cite{KM} that
$\operatorname{Aut}(X_0(48))=\operatorname{Norm}(\Gamma_0(48))/\Gamma_0(N)$.
Now if Claim \ref{AL} where true this group would be isomorphic to
$\Z/2\times \Pi_4$ where $\Pi_n$ is the permutation group of $n$
elements. It is clear that the center of this group is
$\Z/2\times\{1\}$, generated by the Atkin-Lehner involution $w_3$,
but this involution is not the hyperelliptic one.
\end{dem}
The problem of $N=48$ is that $S_4$ does not commute with the
Atkin-Lehner involution $w_3$; thus the direct product
decomposition of Claim \ref{AL} is not possible.

This problem appears also for powers of $3$ one can prove,
\begin{lemma}\label{lema10} Let $N=3^{v_3(N)}\prod_ip_i^{n_i}$ where $p_i$ are different
primes of $\Q$.  Impose that $S_3\in
\operatorname{Norm}(\Gamma_0(N))$. Then $S_3$ commutes with
$w_{p_i^{n_i}}$ if and only if $p_i^{n_i}\equiv 1(modulo\ 3)$.
Therefore if some $p_i^{n_i}\equiv -1(modulo\ 3)$ the Claim
\ref{AL} is not true.
\end{lemma}
\begin{dem} Let us show that $S_{3}$ does not commute with $w_{p_{i}^{n_i}}$
if and only if $p_{i}^{n_i}\equiv-1(mod$ $3$). Observe the
equality
$w_{p_{i}^{n_i}}=\frac{1}{\sqrt{p_{i}^{n_i}}}\left(\begin{array}{cc}
p_{i}^{n_i}k&1\\
Nt&p_{i}^{n_i}\\
\end{array}\right)$:
\begin{center}
$$w_{p_{i}^{n_i}}S_{3}w_{p_{i}^{n_i}}S_{3}^2=$$
$$\frac{1}{p_{i}^{n_i}} \left(
\begin{array}{cc}
(p_{i}^{n_i}k)^2+Nt(1+\frac{p_{i}^{n_i}k}{3})&p_{i}^{n_i}k(\frac{2p_{i}^{n_i}k}{3}+1)+(\frac{p_{i}^{n_i}k}{3}+1)
(\frac{2Nt}{3}+p_{i}^{n_i})\\
Nt(p_{i}^{n_i}k)+Nt(\frac{Nt}{3}+p_{i}^{n_i})&Nt(\frac{2p_{i}^{n_i}k}{3}+1)+p_{i}^{n_i}(\frac{Nt}{3}+p_{i}^{n_i})
(\frac{2Nt}{3}+p_{i}^{n_i})\\
\end{array} \right).$$
\end{center}
For this element to belong to $\Gamma_{0}(N)$ one needs to impose
$\frac{2k^2p_{i}^{n_i}}{3}+\frac{p_{i}^{n_i}k}{3}\in\Z$. Since
$p_{i}^{n_i}\equiv1\ o\ -1(mod\ 3)$ it is needed that
$k\equiv1(mod\ 3)$. Now from $\operatorname{det}(w_{p_i})=1$ we
obtain that $p_{i}^{n_i}k\equiv1(mod\ 3)$; therefore
$p_{i}^{n_i}\equiv1(mod\ 3)$.
\end{dem}

\section{The group structure of $\operatorname{Norm}(\Gamma_0(N))/\Gamma_0(N)$ revisited.}

In this section we correct Claim \ref{AL}. We prove here that the
quotient $$\operatorname{Norm}(\Gamma_0(N))/\Gamma_0(N)$$ is the
product of some groups associated every one of them to the primes
which divide $N$. See for the explicit result theorem
\ref{Barsfi}.

\begin{thm}\label{Bars} Any element $w\in \operatorname{Norm}(\Gamma_0(N))$ has an expression
of the form
$$w=w_{m}\Omega,$$ where $w_m$ is an Atkin-Lehner involution of
$\Gamma_0(N)$ with $(m,6)=1$ and $\Omega$ belongs to the subgroup
generated by $S_{v(N)}$ and the Atkin Lehner involutions
$w_{2^{v_2(N)}}$, $w_{3^{v_3(N)}}$. Moreover for
$\operatorname{gcd}(v(N),2^3)\leq 2$ the group structure for the
subgroup $<S_{v_2(v(N))},w_{2^{v_2(N)}}>$ and $<S_{v_3(v(N))},
w_{3^{v_3(N)}}>$ of $<S_{v(N)},w_{2^{v_2(N)}},w_{3^{v_3(N)}}>$ is
the predicted by Atkin-Lehner at Claim \ref{AL}, but these two
subgroups do not necessary commute withe each other element-wise.
\end{thm}
\begin{dem}
Let us take any element $w$ of the
$\operatorname{Norm}(\Gamma_0(N))$. By Theorem \ref{teo1} we can
express $w$ as follows,
$$w={\sqrt{\delta}}\left(\begin{array}{cc}
r\Delta&\frac{u}{v\delta\Delta}\\
\frac{sN}{v\delta\Delta}&l\Delta\\
\end{array}\right)=\frac{1}{\Delta\sqrt{\delta}}\left(\begin{array}{cc}
r\delta\Delta^2&\frac{u}{v}\\
\frac{sN}{v}&l\delta\Delta^2\\
\end{array}\right)$$

Let us denote by $U=2^{v_2(N)}3^{v_3(N)}$. Write
$\Delta'=\operatorname{gcd}(\Delta,N/U)$ and
$\delta'=\operatorname{gcd}(\delta,N/U)$; then we obtain
$$w_{\delta'{\Delta'}^2}w=\frac{1}{\frac{\Delta}{\Delta'}\sqrt{\delta/\delta'}}\left(\begin{array}{cc}
r'\frac{\delta}{\delta'}\frac{\Delta^2}{\Delta'^2}&\frac{u'}{v(N)}\\
\frac{Nt'}{v(N)}&v'\frac{\delta}{\delta'}\frac{\Delta^2}{\Delta'^2}\\
\end{array}\right)$$
Observe that if $v(N)=1$ we already finish and we reobtain
proposition \ref{propv1}. This is clear if
$\operatorname{gcd}(N,6)=1$; if not, the matrix
$ww_{\delta'\Delta'^2}$ is the Atkin-Lehner involution at
$(\frac{\Delta}{\Delta'})^2\frac{\delta}{\delta'}\in\N$.

Now we need only to check that any matrix of the form
\begin{equation}\label{eq5}\Omega=\frac{1}{\frac{\Delta}{\Delta'}\sqrt{\delta/\delta'}}\left(\begin{array}{cc}
r'\frac{\delta}{\delta'}(\frac{\Delta}{\Delta'})^2&\frac{u'}{v(N)}\\
\frac{Nt'}{v(N)}&v'\frac{\delta}{\delta'}(\frac{\Delta}{\Delta'})^2\\
\end{array}\right)
\end{equation}
 is generated by $S_{v(N)}$ and the Atkin-Lehner involutions at
2 and 3 which are the factors of
$\frac{\delta}{\delta'}(\frac{\Delta}{\Delta'})^2$. To check this
observe that $\Omega=\Omega_2\Omega_3$ with
\begin{scriptsize}
\begin{equation}\label{factor}\Omega_2=\frac{1}{2^{v_2(\frac{\Delta}{\Delta'}\sqrt{\delta/\delta'})}}\left(\begin{array}{cc}
r''2^{v_2(\frac{\delta}{\delta'}(\frac{\Delta}{\Delta'})^2)}&\frac{u''}{2^{v_2(v(N))}}\\
\frac{Nt''}{2^{v_2(v(N))}}&v''2^{v_2(\frac{\delta}{\delta'}(\frac{\Delta}{\Delta'})^2)}\\
\end{array}\right)\end{equation}  $$\Omega_3=\frac{1}{3^{v_3(\frac{\Delta}{\Delta'}\sqrt{\delta/\delta'})}}\left(\begin{array}{cc}
r'''3^{v_3(\frac{\delta}{\delta'}(\frac{\Delta}{\Delta'})^2)}&\frac{u'''}{3^{v_3(v(N))}}\\
\frac{Nt'''}{3^{v_3(v(N))}}&v'''3^{v_3(\frac{\delta}{\delta'}(\frac{\Delta}{\Delta'})^2)}\\
\end{array}\right).$$
\end{scriptsize}
We only consider the case for $\Omega_2$, the case  for the
$\Omega_3$ is similar. We can assume that
$2^{v_2(\frac{\Delta}{\Delta'}\sqrt{\delta/\delta'})}=1$
substituting $\Omega_2$ by $w_{2^{v_2(N)}}\Omega_2$ if necessary.
Thus, we are reduced to a matrix of the form
$\tilde{\Omega_2}=\left(\begin{array}{cc}
r'&\frac{u'}{2^{v_2(v(N))}}\\
\frac{Nt'}{2^{v_2(v(N))}}&v'\\
\end{array}\right)$. Now for some $i$ we can
obtain
$S_{2^{v_2(v(N))}}^i\tilde{\Omega_2}=\left(\begin{array}{cc}
r'&{u'}\\
\frac{Nt'}{2^{v_2(v(N))}}&v'\\
\end{array}\right);$ name this matrix by $\overline{\Omega_2}$.
Then, it is easy to check that $w_{2^{v_2(N)}} S_{2^{v_2(v(N))}}^i
w_{2^{v_2(N)}}\overline{\Omega_2}\in \Gamma_0(N)$ for some $i$.

Similar argument as above are obtained if we multiply $w$ by $w_m$
on the right, i.e. $ww_m$ is also some $\Omega$ as above obtaining
similar conclusion.

Let us see now that the group generated by $S_{v_2(v(N))}$ and the
Atkin-Lehner involutions at 2, and the group generated by
$S_{v_3(v(N))}$ and the Atkin-Lehner involution at 3 have the
structure predicted in Claim \ref{AL} when
$\operatorname{gcd}(v(N),2^3)\leq 2$. We only need to check when
$v(N)$ is a power of 2 or 3 by (\ref{factor}). For $v(N)=1$ the
matrix (\ref{eq5}) is
$w_{\frac{\delta}{\delta'}(\frac{\Delta}{\Delta'})^2}$ (we denote
$w_1:=id$) (we have in this case a much deeper result, see
proposition \ref{propv1}). Take now $v(N)=2$. If
$l=gcd(3,\delta/\delta')$ let $\Omega=w_l\Omega'$; the matrix
$\Omega'$ is as (\ref{eq5}) but with
$\operatorname{gcd}(3,\delta/\delta')=1$, and
$\frac{\delta}{\delta'}\frac{\Delta^2}{\Delta'^2}$ is only a power
of 2. Then $\Omega'\in <S_2,w_{2^{v_2(N)}}>$, let us to precise
the group structure. For $v(N)=2$ we have $v_2(N)=2$ or 3, and we
have already proved the group structure of Claim \cite{AL} in
lemmas \ref{nna:teo},\ref{nne:teo} (we have moreover that Claim
\ref{AL} is true because $S_2$ commutes with the Atkin-Lehner
involutions $w_{p_i^{n_i}}$ if $(p_i,2)=1$, see proposition
\ref{propv2}). Assume now $v(N)=3$. If $l=gcd(2,\delta/\delta')$
and $\Omega=w_l\Omega'$ then $\Omega'$ is as (\ref{eq5}) but with
$gcd(2,\delta/\delta')=1$, and
$\frac{\delta}{\delta'}\frac{\Delta^2}{\Delta'^2}$ is only a power
of 3. Then $\Omega'\in <S_3,w_{3^{v_3(N)}}>$, let us to precise
the group structure. For $v(N)=3$ we have $v_3(N)\geq 2$. Let us
begin with $v_3(N)=2$, then $\Omega'$ is of the form
$$\Omega'=\left(\begin{array}{cc}
r'&\frac{u'}{3}\\
\frac{Nt'}{3}&v'\\
\end{array}\right)=:a(r',u',t',v')$$
(from the formulation of Theorem \ref{teo1} we can consider
$\frac{\Delta}{\Delta'}=1=\frac{\delta}{\delta'}$ because the
factors outside $3$ does not appear if we multiply for a
convenient Atkin-Lehner involution, and for 3 observe that under
our condition $\Delta=1$) and we have
$$a(r',u',t',v')\in\Gamma_{0}(N)\Leftrightarrow t'\equiv u'\equiv0(mod\ 3)$$
$$a(r',u',t',v')w_{9}\in\Gamma_{0}(N)\Leftrightarrow r'\equiv v'\equiv0(mod\ 3)$$
$$a(r',u',t',v')S_{3}\in\Gamma_{0}(N)\Leftrightarrow r'+u'\equiv t'\equiv0(mod\ 3)$$
$$a(r',u',t',v')S_{3}^2\in\Gamma_{0}(N)\Leftrightarrow 2r'+u'\equiv t'\equiv0(mod\ 3)$$
$$a(r',u',t',v')S_{3}w_{9}\in\Gamma_{0}(N)\Leftrightarrow r'\equiv qt'+v'\equiv0(mod\ 3)$$
$$a(r',u',t',v')S_{3}^2w_{9}\in\Gamma_{0}(N)\Leftrightarrow r'\equiv 2qt'+v'\equiv0(mod\ 3)$$
$$a(r',u',t',v')w_{9}S_{3}^2\in\Gamma_{0}(N)\Leftrightarrow r'+u'\equiv v'\equiv0(mod\ 3)$$
$$a(r',u',t',v')w_{9}S_{3}\in\Gamma_{0}(N)\Leftrightarrow r'+2u'\equiv v'\equiv0(mod\ 3)$$
$$a(r',u',t',v')w_{9}S_{3}^2w_{9}\in\Gamma_{0}(N)\Leftrightarrow u'\equiv qt'+v'\equiv0(mod\ 3)$$
$$a(r',u',t',v')S_{3}^2w_{9}S_{3}^2\in\Gamma_{0}(N)\Leftrightarrow u'\equiv 2qt'+v'\equiv0(mod\ 3)$$
$$a(r',u',t',v')S_{3}^2w_{9}S_{3}\in\Gamma_{0}(N)\Leftrightarrow r'+u'\equiv 2t'q+v'\equiv0(mod\ 3)$$
$$a(r',u',t',v')S_{3}w_{9}S_{3}^2\in\Gamma_{0}(N)\Leftrightarrow 2r'+u'\equiv qt'+v'\equiv0(mod\ 3)$$
and these are all the possibilities, proving that the group is
$\{S_3,w_9|S_3^3=w_9^2=(w_9S_3)^3=1\}$ of order 12. Observe that
$S_3$ does not commute with $w_2$ (see for example lemma
\ref{nna:teo}).

Suppose now that $v_3(N)\geq 3$. We distinguish the cases $v_3(N)$
odd and $v_3(N)$ even. Suppose $v_3(N)$ is even, then
$\frac{\delta}{\delta'}=1$ and $\Omega'$ has the following form
$$\frac{1}{\frac{\Delta}{\Delta'}}\left(\begin{array}{cc}
r'(\frac{\Delta}{\Delta'})^2&\frac{u'}{3}\\
\frac{Nt'}{3}&v'(\frac{\Delta}{\Delta'})^2\\
\end{array}\right)$$
with $\alpha:=\Delta/\Delta'$ dividing $3^{[v_3(N)/2]-1}$. Since
this last matrix has determinant 1 we see that $\alpha$ satisfies
$\operatorname{gcd}(\alpha,N/(3^2\alpha^2))=1$; thus $\alpha=1$ or
$\alpha=3^{[v_3(N)/2]-1}$. Write
$a(r',u',t',v')=\left(\begin{array}{cc}
r'&\frac{u'}{3}\\
\frac{Nt'}{3}&v'\\
\end{array}\right)$ when we take
$\alpha=1$ and $b(r',u',t',v')=\left(\begin{array}{cc}
r'(3^{[v_3(N)/2]-1})&\frac{u'}{3^{[v_3(N)/2]}}\\
\frac{Nt'}{3^{[v_3(N)/2]}}&v'(3^{[v_3(N)/2]-1})\\
\end{array}\right)$ when
$\alpha=3^{[v_3(N)/2]-1}$. It is easy to check that
$b(r',u',t',v')=w_{3^{v_3(N)}}a(r',u',t',v')$ and that the group
structure is the predicted in a similar way as the one done above
for $v(N)=2$. Suppose now that $v_3(N)$ is odd, then
$\frac{\delta}{\delta'}$ is 1 or 3 and $\frac{\Delta}{\Delta'}$
divides $3^{[v_3(N)/2]-1}$. Now from $\operatorname{det}()=1$ we
obtain that the only possibilities are
$\frac{\delta}{\delta'}=1=\frac{\Delta}{\Delta'}$ name the
matrices for this case following equation \ref{eq5} by
$a(r',u',t',v')$, and the other possibility is
$\frac{\delta}{\delta'}=3$ and
$\frac{\Delta}{\Delta'}=3^{[v_3(N)/2]-1}$, write the matrices for
this case following equation \ref{eq5} by $c(r',u',t',v')$. It is
also easy to check that
$c(r',u',t',v')=w_{3^{v_3(N)}}a(r'',u'',t'',v'')$, and that the
group structure is the predicted.

\end{dem}

\begin{cor}\label{propv3}
Let $N=3^{v_3(N)}\prod_{i}p_{i}^{n_i}$, with $p_{i}$ different
primes such that $\operatorname{gcd}(p_{i},6)=1$. Suppose that
$v(N)=3$ and $p_{i}^{n_i}\equiv1(mod\ 3)$ for all $i$. Then Claim
\ref{AL} is true.
\end{cor}
\begin{dem} From the proof of the above theorem \ref{Bars} for $v(N)=3$
with $v_3(N)\geq 2$, lemma \ref{lema10}, and that the general
observation that the Atkin-Lehner involutions commute one with
each other we obtain that the direct product decomposition of
Claim \ref{AL} is true obtaining the result.
\end{dem}

Now we shows the corrections to Claim \ref{AL} for $v(N)=4$ and
$v(N)=8$, about the group structure of the subgroup of
$\operatorname{Norm}(\Gamma_0(N))/\Gamma_0(N)$ generated for
$S_{2^k}$ and the Atkin-Lehner involution at prime 2.

\begin{prop} Suppose $v(N)=4$, observe that in this situation $v_2(N)=4,$ or $5$. Then the group structure of the
subgroup  $<w_{2^{v_2(N)}},S_4>$ of
$\operatorname{Norm}(\Gamma_0(N))/\Gamma_0(N)$ is given by the
relations:
\begin{enumerate}
\item For $v_2(N)=4$ we have $S_4^4=w_{16}^{2}=(w_{16}S_4)^3=1$.
\item For $v_2(N)=5$ we have $S_4^4=w_{32}^{2}=(w_{32}S_4)^4=1$.
\end{enumerate}
\end{prop}
\begin{dem} It is a straightforward computation. Observe that for
$v_2(N)=4$ the statement coincides with Claim \ref{AL} but not for
$v_2(N)=5$, where one checks that $S_4$ does not commute with
$w_{32}S_4w_{32}$.
\end{dem}
\begin{prop} Suppose $v(N)=8$ and $v_2(N)$ even (this is the case (3)(c) in Claim \ref{AL}).
Then the group
$<w_{2^{v_2(N)}},S_8>\subseteq\operatorname{Norm}(\Gamma_0(N))/\Gamma_0(N)$
satisfies the following relations: $S_8^8=w_{2^{v_2(N)}}^2=1$, and
\begin{enumerate}
\item for $v_2(N)=6$ we have $(w_{64}S_8)^3=1$, \item for
$v_2(N)\geq 8$ we do not have the relation
$(w_{2^{v_2(N)}}S_8)^3=1$, \item for $v_2(N)\geq 10$ we have the
relation: $S_8$ commutes with $w_{2^{v_2(N)}}S_8 w_{2^{v_2(N)}}$,
\item for $v_2(N)=6$ or $8$ we do not have the relation: $S_8$
commutes with the element $w_{2^{v_2(N)}}S_8 w_{2^{v_2(N)}}$.
\item For $v_2(N)=8$ we have the relation:
$w_{256}S_8w_{256}S_8w_{256}S_8^3w_{256}S_8^3=1$.
\end{enumerate}
\end{prop}
\begin{dem} Straightforward.
\end{dem}
\begin{prop} Suppose $v(N)=8$ and $v_2(N)$ odd (this is the case (3)(d) in Claim \ref{AL}).
Then the group
$<w_{2^{v_2(N)}},S_8>\subseteq\operatorname{Norm}(\Gamma_0(N))/\Gamma_0(N)$
satisfies the following relations: $S_8^8=w_{2^{v_2(N)}}^2=1$, and
\begin{enumerate}
\item for $v_2(N)=7$ $(w_{128}S_8)^4=1$, \item for $v_2(N)\geq 9$
we do not have the relation $(w_{2^{v_2(N)}}S_8)^4=1$, \item for
$v_2(N)\geq 9$ we have the Atkin-Lehner relation:
%\begin{enumerate}
%\item $S_8^7$ commutes with $w_{2^{v_2(N)}}S_8^5 w_{2^{v_2(N)}}$
 $S_8$ commutes with $w_{2^{v_2(N)}}S_8 w_{2^{v_2(N)}}$,
 \item for $v_2(N)=7$ we do not have that $S_8$ commutes with
 $w_{128}S_8w_{128}$.
%\end{enumerate}
\end{enumerate}
\end{prop}
\begin{dem} Straightforward.
\end{dem}
Let us finally write the revisited results concerning Claim
\ref{AL} that we prove;
\begin{thm}\label{Barsfi} The quotient $\operatorname{Norm}(\Gamma_{0}(N))/\Gamma_0(N)$ is
a product of the following groups:
\begin{enumerate}
\item \{$ w_{q^{\upsilon_q(N)}}$\} for every prime $q$, $q\ge5$
$q\mid N$. \item
\begin{enumerate}
      \item If $\upsilon_{3}(N)=0$, \{$1$\}
      \item If $\upsilon_{3}(N)=1$, \{$w_{3}$\}
      \item If $\upsilon_{3}(N)=2$, \{$w_{9},S_{3}$\}; satisfying $w_{9}^2=S_{3}^3=(w_{9}S_{3})^3=1$ (factor of order 12)
      \item If $\upsilon_{3}(N)\ge3$; \{$w_{3^{\upsilon_{3}(N)}},S_{3}$\}; where $w_{3^{\upsilon_{3}(N)}}^2=S_{3}^3=1$ and $w_{3^{\upsilon_{3}(N)}}S_{3}w_{3^{\upsilon_{3}(N)}}$
      commute
      with $S_{3}$ (factor group with 18 elements)
      \end{enumerate}
\item Let be $\lambda=\upsilon_{2}(N)$ and
$\mu=\operatorname{min}(3,[\frac{\lambda}{2}])$ and denote by
$\upsilon''=2^{\mu}$ the we have:
      \begin{enumerate}
      \item If $\lambda=0$ ; \{$1$\}
      \item If $\lambda=1$; \{$w_{2}$\}
      \item If $\lambda=2\mu$ and $2\leq\lambda\leq 6$; \{$w_{2^{\upsilon_{2}(N)}},S_{\upsilon''}$\}
      with the relations
      $w_{2^{\upsilon_{2}(N)}}^2=S_{\upsilon''}^{\upsilon''}=(w_{2^{\upsilon_{2}(N)}}S_{\upsilon''})^3=1$,
      where they have orders 6,24, and 96 for $\upsilon=2,4,8$ respectively.
      \item If $\lambda> 2\mu$ and $2\leq\lambda\leq 7$;
      \{ $w_{2^{\upsilon_{2}(N)}},S_{\upsilon''}$\};
      $w_{2^{\upsilon_{2}(N)}}^2=S_{\upsilon''}^{\upsilon''}=1$.
       Moreover, $(w_{2^{\upsilon_{2}(N)}}S_{\upsilon''})^4=1$.
       \item[($\tilde{c}$),($\tilde{d}$)] If $\lambda\geq 9$; \{$w_{2^{\upsilon_{2}(N)}},S_{8}$\}
      with the relations
      $w_{2^{\upsilon_{2}(N)}}^2=S_{8}^{8}=1$ and $S_8$
       commutes with $w_{2^{v_2(N)}}S_8w_{2^{v_2(N)}}$.
       \item[($\hat{c}$)] If $\lambda=8$; \{$w_{2^{\upsilon_{2}(N)}},S_{8}$\}
      with the relations
      $w_{2^{\upsilon_{2}(N)}}^2=S_{8}^{8}=1$
      and $w_{256}S_8w_{256}S_8w_{256}S_8^3w_{256}S_8^3=1$.
      \end{enumerate}
\end{enumerate}
\end{thm}
\begin{obs2} One needs to warn that for the situation $v(N)=8$
possible the relations does not define totally the factor group,
but it is a computation more.
\end{obs2}
\begin{obs2} The product between the different groups appearing in
theorem \ref{Barsfi} is easily computable. Effectively, we know
that the Atkin-Lehner involutions commute, and $S_{2^{v_2(v(N))}}$
commutes with $S_{3^{v_3(v(N))}}$. Moreover $S_2$ commutes with
any element from lemma \ref{lema2}. Consider $w_{p^n}$ an
Atkin-Lehner involution for $X_0(N)$ with $p$ a prime. One obtains
the following results by using the same arguments appearing in the
proof of lemma \ref{lema10};
\begin{enumerate}
\item let $p$ be coprime with $3$ and $3|v(N)$. $S_3$ commutes
with $w_{p^n}$ if and only if $p^n\equiv 1(modulo\ 3)$. If
$p^n\equiv-1(modulo\ 3)$ then $w_{p^n}S_3=S_3^2w_{p^n}$. \item Let
$p$ be coprime with $2$ and $4|v(N)$. $S_4$ commutes with
$w_{p^n}$ if and only if $p^n\equiv 1(modulo\ 4)$. If
$p^n\equiv-1(modulo\ 4)$ then $w_{p^n}S_4=S_4^3w_{p^n}$. \item Let
$p$ be coprime with $2$ and $8|v(N)$. Then,
$w_{p^n}S_8=S_8^kw_{p^n}$ if $p^n\equiv k(modulo\ 8)$, in
particular $S_8$ commutes with $w_{p^n}$ if and only if $p^n\equiv
1(modulo\ 8)$.

\end{enumerate}
\end{obs2}

\vspace{1cm}

Francesc Bars Cortina, Depart. Matem\`atiques, Universitat
Aut\`onoma de Barcelona, 08193 Bellaterra. Catalonia. Spain.
E-mail:
francesc@mat.uab.cat \\

\end{document}